\documentclass[12pt,notitlepage,twoside]{article}
\usepackage{latexsym} 
\usepackage{amsmath} 
\usepackage{amsfonts} 
\usepackage{amssymb} 
\renewcommand{\a }{\alpha }

\newcommand{\D }{\Delta }
\newcommand{\e }{\varepsilon }

\renewcommand{\l }{\lambda }

\newcommand{\n }{\nabla }

\newcommand{\hes }{{\rm Hess}\,}

\newcommand{\wtilde }{\widetilde} 
 
\newcommand{\be}{\begin{equation}} 
\newcommand{\ee}{\end{equation}} 
\newenvironment{pf}{\noindent{\bf Proof.}\enspace}{%\rule{2mm}{2mm}
\hfill$\Box$\medskip} 
\newenvironment{pfn}[1]{\noindent{\bf Proof of {#1}\enspace}}{%\rule{2mm}{2mm}
\hfill$\Box$\medskip} 
\newcommand{\R}{\mathbb{R}}

\newtheorem{thm}{Theorem}[section] 
\newtheorem{pro}[thm]{Proposition}
\newtheorem{lem}[thm]{Lemma}
\newtheorem{rem}[thm]{Remark}

\numberwithin{equation}{section}

\author{{\sc  Veronica Felli}
 \thanks{Supported by MIUR, national project ``Variational Methods and
   Nonlinear Differential Equations''.}}

\title { \Large \textbf{A note on the existence of $H$-bubbles\\
 via perturbation methods}} 
 
\begin{document}

\date{S.I.S.S.A. - Via Beirut 2-4\\ 
34014 Trieste, Italy\\
e-mail: \texttt{felli@sissa.it}}

\maketitle

{\footnotesize 
\begin{abstract} 
 
\noindent 
We study the problem of existence of surfaces in $\R^3$ parametrized
on the sphere~${\mathbb S}^2$ with prescribed mean curvature $H$ in
the perturbative case, i.e. for $H=H_0+\e
H_1$, where $H_0$ is a nonzero constant, $H_1$ is a $C^2$ function and
$\e$ is a small perturbation parameter.

\bigskip\bigskip



\end{abstract}

}

\section{Introduction}
\mbox{}

\noindent In this paper we are interested in the existence of $H$-bubbles, namely of ${\mathbb S}^2$-type parametric surfaces in $ \R^3$ with prescribed mean curvature $H$. This geometrical problem is motivated by some models describing capillarity phenomena and has the following analytical formulation: given a function $H\in C^1(\R^3)$, find a smooth nonconstant function $\omega: \R^2\to\R^3$ which is conformal as a map on ${\mathbb S}^2$ and solves the problem
\begin{equation}\tag{$P_H$}
\begin{cases}
\Delta \omega=2 H(\omega)\,\omega_x\wedge\omega_y,&\text{in }\R^2,\\[5pt]
{\displaystyle \int_{\R^2}|\n \omega|^2<+\infty},
\end{cases}
\end{equation}
where $\omega_x=\big(\frac{\partial\omega_1}{\partial
  x},\frac{\partial\omega_2}{\partial
  x},\frac{\partial\omega_3}{\partial x}\big)$,
$\omega_y=\big(\frac{\partial\omega_1}{\partial
  y},\frac{\partial\omega_2}{\partial
  y},\frac{\partial\omega_3}{\partial y}\big)$,
$\D\omega=\omega_{xx}+\omega_{yy}$, $\n \omega=(\omega_x,\omega_y)$,
and $\wedge$ denotes the exterior product in $\R^3$.\par
Brezis and Coron \cite{BC} proved that for constant nonzero mean
curvature $H(u)\equiv H_0$ the only nonconstant solutions are spheres
of radius $|H_0|^{-1}$. \par
 While the Plateau and the Dirichlet problems  has been largely studied
 both for $H$ constant and for $H$
 nonconstant (see \cite{bc, BC, Hi, Jak1, Jak2, steffen, stru,
   struwe}), problem $(P_H)$ in the case of  nonconstant $H$ has been
 investigated 
 only recently, see \cite{CM1, CM2, CM}. In \cite{CM1} Caldiroli and
Musina proved the existence of $H$-bubbles with minimal energy under
the assumptions that $H\in C^1(\R^3)$ satisfies 
\begin{align*}
(i)\quad&\sup_{u\in\R^3}|\n H(u+\xi)\cdot u \,u|<1\quad\text{for some }\xi\in\R^3,\\
(ii)\quad&H(u)\to H_{\infty}\quad\text{as } |u|\to\infty\quad\text{for some }H_{\infty}\in\R,\\
(iii)\quad&c_H=\inf_{\substack{u\in C^1_c(\R^2,\R^3)\\u\neq 0}}\sup_
  {s>0}{\mathcal E}_H(su)<\frac{4\pi}{3H_{\infty}^2}
\end{align*}
where ${\mathcal E}_H(u)=\frac12 \int_{\R^2}|\n
u|^2+2\int_{\R^2}Q(u)\cdot u_x\wedge u_y$ and $Q:\R^3\to\R^3$ is any vector field such that ${\rm div\,}Q=H$.\par
The perturbative method introduced by Ambrosetti and Badiale
\cite{AB1,AB2} was used in~\cite{CM} to treat the case in which $H$ is
a small perturbation of a constant, namely 
\[
H(u)=H_{\e}(u)=H_0+\e H_1(u),
\]
 where $H_0\in\R\setminus\{0\}$, $H_1\in C^2(\R^3)$, and $\e$ is a
 small real parameter. This method allows to find critical points of a
 functional $f_{\e} $ of the type $f_{\e}(u)=f_0(u)-\e G(u)$ in a
 Banach space by studying a finite dimensional problem. More
 precisely, if the unperturbed functional $f_0$ has a finite
 dimensional manifold of critical points $Z$ which satisfies a
 nondegeneracy condition, it is possible to prove, for $|\e|$
 sufficiently small, the existence of a smooth function
 $\eta_{\e}(z):Z\to(T_zZ)^{\perp}$ such that any critical point $\bar
 z\in Z$ of the function 
\[
\Phi_{\e}:Z\to\R,\quad \Phi_{\e}(z)=f_{\e}\big(z+\eta_{\e}(z)\big)
\]
gives rise to a critical point $u_{\e}=\bar z+\eta_{\e}(\bar z)$ of $f_{\e}$, i.e. the perturbed manifold $Z_{\e}:=\{z+\eta_{\e}(z):\ z\in Z\}$ is a natural constraint for $f_{\e}$. Furthermore $\Phi_{\e}$ can be expanded as 
\begin{equation}\label{eq:exp}
\Phi_{\e}(z)=b-\e \Gamma(z)+o(\e)\quad\text{as }\e\to 0
\end{equation}
where $b=f_0(z)$ and $\Gamma$ is the Melnikov function defined as the restriction of the perturbation $G$ on $Z$, namely $\Gamma=G\big|_Z$. For problem $(P_{H_{\e}})$, $\Gamma$
 is given by 
\[
\Gamma:\R^3\to\R,\quad \Gamma(p)=\int_{|p-q|<\frac1{|H_0|}}H_1(q)\,dq.
\]
In \cite{CM} Caldiroli and Musina studied the functional $\Gamma$
giving some existence results in the perturbative setting for problem
$(P_{H_{\e}})$. They prove that for $|\e|$ small there exists a smooth
$H_{\e}$-bubble if one of the following conditions
holds
\begin{align*}
1)\quad& H_1 \text{ has a 
nondegenerate stationary point and } |H_0| \text{ is large;}\\
2)\quad& \max_{p\in\partial K}H_1(p)<\max_{p\in
  K}H_1(p)\quad\text{or}\min_{p\in\partial K}H_1(p)>\min_{p\in
  K}H_1(p)\\
&\text{for some nonempty compact set } K\subset\R^3 \text{
  and }|H_0| \text{ is large;}\\
3)\quad& H_1\in L^r(\R^3)\text{ for some } r\in [1,2].
\end{align*}
They prove that critical points of $\Gamma$ give rise to solutions to
$(P_{H_{\e}})$ for $\e$ sufficiently small. Precisely the assumption
that $H_0$ is large required in cases 
1) and 2) ensures that if $H_1$ is not constant then $\Gamma$ is not
identically constant. If we let this assumption drop, it may happen
that $\Gamma$ is constant even if $H_1$ is not. This fact produces
some loss of information because the first 
order expansion \eqref{eq:exp} is not sufficient to deduce the existence
of critical points of $\Phi_{\e}$ from the existence of critical
points of $\Gamma$. Instead of studying $\Gamma$ we perform a direct
study of $\Phi_{\e}$ which allows us to prove some new results. In
the first one, we assume that $H_1$ 
vanishes at $\infty$ and has bounded gradient, and prove the existence
of a solution without branch points. Let us recall that a branch point
for a solution $\omega$ to $(P_H)$ is a point where $\n \omega=0$,
i.e. a point where the surface parametrized by $\omega$ fails to be
immersed.  
\begin{thm}\label{t:1}
Let $H_0\in\R\setminus\{0\}$, $H_1\in C^2(\R^3)$ such that
\begin{align*}
(H1)\quad &\lim_{|p|\to\infty}H_1(p)=0;\\
(H2)\quad &\n H_1\in L^{\infty}(\R^3,\R^3).
\end{align*}
Let $H_{\e}=H_0+\e H_1$. Then for $|\e|$ sufficiently small there
exists a smooth $H_{\e}$-bubble without branch points.
\end{thm}
With respect to case 1) of \cite{CM} we require neither nondegeneracy of
critical points of $H_1$ nor largeness of $H_0$. With respect to case 2)
we do not assume that $H_0$ is large; on the other hand our assumption
$(H1)$ implies 2). Moreover we do not assume any integrability condition of
type 3).  With respect to the result proved in \cite{CM1}, we have the
same kind of behavior of $H_1$ at $\infty$ (see $(ii)$ and assumption
$(H1)$) but we do not need any assumption of type $(iii)$; on the
other hand in \cite{CM1} it is not required that the prescribed curvature is a
small perturbation of a constant. 

The following results give some conditions on $H_1$ in
order to have two or three solutions.
\begin{thm}\label{t:2}
Let $H_0\in\R\setminus\{0\}$, $H_1\in C^2(\R^3)$ such that $(H1)$,
$(H2)$, 
\begin{align*}
(H3)\quad &\hes H_1(p) \text{ is positive definite for any } p\in
B_{1/|H_0|}(0),\\
(H4)\quad &H_1(p)>0 \text{ in } B_{1/|H_0|}(0),
\end{align*}
hold. Then for $|\e|$ sufficiently small there exist at least three
smooth $H_{\e}$-bubbles without branch points. 
\end{thm}
\begin{rem}\label{r:2}
If we assume $(H1)$, $(H2)$, and, instead of $(H3)-(H4)$, that
$H_1(0)>0$ and $\hes
H_1(0)$ is positive definite, then we can
prove that for $|H_0|$ sufficiently large and $|\e|$ sufficiently
small there exist at least three smooth $H_{\e}$-bubbles without branch points.
\end{rem}
\begin{thm}\label{t:3}
Let $H_0\in\R\setminus\{0\}$, $H_1\in C^2(\R^3)$ such that $(H1)$ and
$(H2)$ hold. Assume that there exist $p_1,p_2\in\R^3$ such that
\[
(H5)\quad \int_{B(p_1,1/|H_0|)}H_1(\xi)\,d\xi>0\quad\text{and}\quad
\int_{B(p_2,1/|H_0|)}H_1(\xi)\,d\xi<0.
\]
Then for $|\e|$ sufficiently small there exist at least two smooth
$H_{\e}$-bubbles without branch points.
\end{thm}
\begin{rem}\label{r:3}
If we assume $(H1)$, $(H2)$, and, instead of $(H5)$, that there exist
$p_1,p_2\in\R^3$ such that $H_1(p_1)>0$ and $H_1(p_2)<0$, then we can
prove that for $|H_0|$ sufficiently large and $|\e|$ sufficiently
small there exist at least two smooth $H_{\e}$-bubbles without branch points.
\end{rem}
The present paper is organized as follows. In Section 2 we introduce
some notation and recall some known facts whereas Section 3 is devoted
to the proof of Theorems \ref{t:1}, \ref{t:2}, and~\ref{t:3}.

\medskip\noindent
{\bf Acknowledgments.}\quad  The author wishes to thank Professor
A. Ambrosetti and Professor R. Musina for many helpful suggestions.

\section{Notation and known facts}

\mbox{}

\noindent In the sequel we will take $H_0=1$; this is not restrictive since we can do the change $H_1(u)=H_0\wtilde H_1(H_0u)$. Hence we will always write
\[
H_{\e}(u)=1+\e H(u),
\]
where $H\in C^2(\R^3)$. Let us denote by $\omega$ the function $\omega:\R^2\to{\mathbb S}^2$ defined as
\[
\omega(x,y)=\big(\mu(x,y)x,\mu(x,y)y,1-\mu(x,y)\big)\quad\text{where}\quad  \mu(x,y)=\frac2{1+x^2+y^2},\ (x,y)\in \R^2.
\]
Note that $\omega$ is a conformal parametrization of the unit sphere and solves
\begin{equation}\label{eq:omega}
\begin{cases}
\Delta\omega=2\,\omega_x\wedge\omega_y&\text{on }\R^2\\[5pt]
{\displaystyle \int_{\R^2}|\n \omega|^2<+\infty}.
\end{cases}
\end{equation}
Problem \eqref{eq:omega} has in fact a family of solutions of the form $\omega\circ\phi+p$ where $p\in\R^3$ and $\phi$ is any conformal diffeomorphism of $\R^2\cup\{\infty\}$. For $s\in(1,+\infty)$, we will set $L^s:=L^s({\mathbb S}^2,\R^3)$, where any map $v\in L^s$ is identified with the map $v\circ \omega:\R^2\to\R^3$ which satisfies
\[
\int_{\R^2}|v\circ\omega|^s\mu^2=\int_{{\mathbb S}^2}|v|^s.
\]
We will use the same notation for $v$ and $v\circ\omega$. By $W^{1,s}$ we denote the Sobolev space $W^{1,s}({\mathbb S}^2,\R^3)$ endowed (according to the above identification) with the norm
\[
\|v\|_{W^{1,s}}=\bigg[\int_{\R^2}|\n v|^s\mu^{2-s}\bigg]^{1/s}+\bigg[\int_{\R^2}|v|^s\mu^{2}\bigg]^{1/s}.
\]
If $s'$ is the conjugate exponent of $s$, i.e. $\frac1s+\frac 1{s'}=1$, the duality product between $W^{1,s}$ and $W^{1,s'}$ is given by
\[
\langle v,\varphi\rangle=\int_{\R^2}\n v\cdot\n \varphi+\int_{\R^2} v\cdot\varphi\,\mu^2\quad\text{for any }v\in W^{1,s}\text{ and } \varphi\in W^{1,s'}.
\]
Let $Q$ be any smooth vector field on $\R^3$ such that ${\rm div\,}Q=H$. The energy functional associated to problem 
\begin{equation}\tag{$P_{\e}$}
\begin{cases}
\Delta u=2 \big(1+\e H(u)\big)\,u_x\wedge   u_y,&\text{in }\R^2,\\[5pt]
{\displaystyle \int_{\R^2}|\n u|^2<+\infty},
\end{cases}
\end{equation}
is given by 
\[
{\mathcal E}_{\e}(u)=\frac12 \int_{\R^2}|\n u|^2+2{\mathcal V}_1(u)+2\e{\mathcal V}_H(u),\quad u\in W^{1,3},
\]
where 
\[
{\mathcal V}_H(u)=\int_{\R^2}Q(u)\cdot u_x\wedge u_y
\]
has the meaning of an algebraic volume enclosed by the surface
parametrized by $u$ with weight $H$ (it is independent of the choice of $Q$); in particular
\[
{\mathcal V}_1(u)=\frac13\int_{\R^2}u\cdot u_x\wedge u_y.
\]
In \cite{CM}, Caldiroli and Musina studied some regularity properties
of ${\mathcal V}_H$ on the space $W^{1,3}$. In particular they proved
the following properties. 
\begin{itemize}
\item[a)] For $H\in C^1(\R^3)$, the functional ${\mathcal V}_H$ is of
  class $C^1$ on $W^{1,3}$ and the Fr\'echet differential of
  ${\mathcal V}_H$ at $u\in W^{1,3}$ is given by 
\begin{equation}\label{eq:frechet}
d{\mathcal V}_H(u)\varphi=\int_{\R^2}H(u)\,\varphi\cdot u_x\wedge
u_y\quad\text{for any }\varphi\in W^{1,3} 
\end{equation}
and admits a unique continuous and linear extension on $W^{1,3/2}$
defined by \eqref{eq:frechet}. Moreover for every $u\in W^{1,3}$ there
exists ${\mathcal V}_H'(u)\in W^{1,3}$ such that  
\begin{equation}\label{eq:frechet1}
\langle{\mathcal V}_H'(u),\varphi\rangle=\int_{\R^2}H(u)\,\varphi\cdot
u_x\wedge u_y\quad\text{for any }\varphi\in W^{1,3/2}. 
\end{equation}
\item[b)] For $H\in C^2(\R^3)$, the map ${\mathcal V}_H':W^{1,3}\to W^{1,3}$ is of class $C^1$ and 
\begin{align}\label{eq:frechet2}
\langle {\mathcal
  V}_H''(u)\cdot\eta,\varphi\rangle=\int_{\R^2}H(u)\,
\varphi\cdot&(\eta_x\wedge u_y+u_x\wedge \eta_y)+\int_{\R^2}(\n
H(u)\cdot\eta)\varphi\cdot(u_x\wedge u_y)\notag\\   
&\text{for any }u,\eta\in W^{1,3}\text{ and }\varphi\in W^{1,3/2}. 
\end{align}
\end{itemize}
Hence for all $u\in W^{1,3}$, ${\mathcal E}_{\e}'(u)\in W^{1,3}$ and for any $\varphi\in W^{1,3/2}$
\[
\langle  {\mathcal E}_{\e}'(u),\varphi\rangle=\int_{\R^2}\n
u\cdot\n\varphi+2\int_{\R^2}\varphi\cdot u_x\wedge
u_y+2\e\int_{\R^2}H(u)\,\varphi\cdot u_x\wedge u_y. 
\]
As remarked in \cite{CM}, critical points of ${\mathcal E}_{\e}$ in $W^{1,3}$ give rise to bounded weak solutions to $(P_{\e})$ and hence by the regularity theory for $H$-systems (see \cite{gruter}) to classical conformal solutions which are $C^{3,\alpha}$ as maps on ${\mathbb S}^2$.\par
The unperturbed problem, i.e. $(P_{\e})$ for $\e=0$, has a 9-dimensional manifold of solutions given by
\[
Z=\{R\omega\circ L_{\l,\xi}+p:\ R\in SO(3),\ \l>0,\ \xi\in \R^2,\ p\in\R^3\}
\]
where $L_{\l,\xi}z=\l(z-\xi)$ (see \cite{isobe}). In \cite{isobe} the
nondegeneracy condition $T_uZ=\ker {\mathcal E}''_0(u)$ for any $u\in
Z$ (where $T_uZ$ denotes the tangent space of $Z$ at $u$) is
proved (see also \cite{ChMa}).\par
As observed in \cite{CM}, in performing the finite dimensional reduction, the dependence on the 6-dimensional conformal group can be neglected since any $H$-system is conformally invariant. Hence we look for critical points of ${\mathcal E}_{\e}$ constrained on a three-dimensional manifold $Z_{\e}$ just depending on the translation variable $p\in\R^3$.

\section{Proof of Theorem \ref{t:1}}
\mbox{}

\noindent We start by constructing a perturbed manifold which is a natural constraint for ${\mathcal E}_{\e}$. 
\begin{lem}\label{l:eta}
Assume $H\in C^2(\R^3)\cap L^{\infty}(\R^3)$ and $\n H\in
L^{\infty}(\R^3,\R^3)$. Then there exist $\e_0>0$, $C_1>0$, and a
$C^1$ map $\eta:(-\e_0,\e_0)\times\R^3\to W^{1,3}$ such that for any
$p\in \R^3$ and 
$\e\in(-\e_0,\e_0)$
\begin{align}
&{\mathcal E}_{\e}'\big(\omega+p+\eta(\e,p)\big)\in T_{\omega}Z,\label{eq:eta1}\\
&\eta(\e,p)\in (T_{\omega}Z)^{\perp},\label{eq:eta2}\\
&\int_{{\mathbb S}^2}\eta(\e,p)=0,\label{eq:eta3}\\
&\|\eta(\e,p)\|_{W^{1,3}}\leq C_1|\e|.\label{eq:eta4}
\end{align}
Moreover if we assume that the limit of $H$ at $\infty$ exists and 
\begin{equation}\label{eq:limitH}
\lim_{|p|\to\infty}H(p)=0
\end{equation}
we have that $\eta(\e,p)$ converges to $0$ in $W^{1,3}$ as $|p|\to\infty$ uniformly with respect to $|\e|< \e_0$.
\end{lem}
\begin{pf}
Let us define the map 
\begin{align*}
&F=(F_1,F_2):\R\times\R^3\times W^{1,3}\times\R^6\times\R^3\to W^{1,3}\times\R^6\times\R^3\\
&\langle F_1(\e,p,\eta,\l,\a), \varphi\rangle=\langle {\mathcal E}_{\e}'(\omega+p+\eta),\varphi\rangle-\sum_{i=1}^6\l_i\int_{\R^2}\n \varphi\cdot\n\tau_i+\a\cdot\int_{{\mathbb S}^2}\varphi,\quad\forall \varphi\in W^{1,3/2}\\
& F_2(\e,p,\eta,\l,\a)=\bigg(\int_{\R^2}\n \eta\cdot \n\tau_1,\dots,\int_{\R^2}\n \eta\cdot \n\tau_6,\int_{{\mathbb S}^2}\eta\bigg)
\end{align*}
where $\tau_1,\dots,\tau_6$ are chosen in $T_{\omega}Z$ such that
\[
\int_{\R^2}\n \tau_i\cdot \n\tau_j=\delta_{ij}\quad\text{and}\quad \int_{{\mathbb S}^2}\tau_i=0\quad i,j=1,\dots,6
\]
so that $T_{\omega}Z$ is spanned by $\tau_1,\dots,\tau_6,e_1,e_2,e_3$. It has been proved by Caldiroli and Musina~\cite{CM} that $F$ is of class $C^1$ and that the linear continuous operator 
\begin{align*}
&{\mathcal L}: W^{1,3}\times\R^6\times\R^3\to W^{1,3}\times\R^6\times\R^3\\
&{\mathcal L}=\frac{\partial F}{\partial(\eta,\l,\a)}(0,p,0,0,0)
\end{align*}
i.e.
\begin{align*}
\langle {\mathcal L}_1(v,\mu,\beta),\varphi\rangle&=\langle {\mathcal E}_0''(\omega)\cdot v,\varphi\rangle-\sum_{i=1}^6\mu_i\int_{\R^2}\n \varphi\cdot\tau_i-\beta\int_{{\mathbb S}^2}\varphi\quad \forall\,\varphi\in W^{1,3/2}\\
{\mathcal L}_2(v,\mu,\beta)&=\bigg(\int_{\R^2}\n v\cdot \n\tau_1,\dots,\int_{\R^2}\n v\cdot \n\tau_6,\int_{{\mathbb S}^2}v\bigg)
\end{align*}
is invertible. Caldiroli and Musina applied the Implicit Function Theorem to solve the equation $F(\e,p,\eta,\l,\a)=0$ locally with respect to the variables $\e,p$, thus finding a $C^1$-function $\eta$ on a neighborhood $(-\e_0,\e_0)\times B_R\subset \R\times\R^3$ satisfying \eqref{eq:eta1}, \eqref{eq:eta2}, and \eqref{eq:eta3}. We will use instead the Contraction Mapping Theorem, which allows to prove the existence of such a function $\eta$ globally on $\R^3$, thanks to the fact that the operator ${\mathcal L}$ does not depend on $p$ and hence it is invertible uniformly with respect to $p\in\R^3$. \par
We have that  $F(\e,p,\eta,\l,\a)=0$ if and only if $(\eta,\l,\a)$ is a fixed point of the map $T_{\e,p}$ defined as
\[
T_{\e,p}(\eta,\l,\a)=-{\mathcal L}^{-1} F(\e,p,\eta,\l,\a)+(\eta,\l,\a).
\]
To prove the existence of $\eta$ satisfying \eqref{eq:eta1}, \eqref{eq:eta2}, and \eqref{eq:eta3}, it is enough to prove that $T_{\e,p}$ is a contraction in some ball $B_{\rho}(0)$ with $\rho=\rho({\e})>0$ independent of $p$, whereas the regularity of $\eta(\e,p)$ follows from the Implicit Function Theorem. \par
We have that if $\|\eta\|_{W^{1,3}}\leq\rho$
\begin{align}\label{eq:1}
\|T_{\e,p}&(\eta,\l,\a)\|_{W^{1,3}\times\R^6\times\R^3}\notag\\
&\leq C_2\|F(\e,p,\eta,\l,\a)-{\mathcal
  L}(\eta,\l,\a)\|_{W^{1,3}\times\R^6\times\R^3}\leq C_2\|{\mathcal
  E}_{\e}'(\omega+p+\eta)-{\mathcal
  E}_0''(\omega)\eta\|_{W^{1,3}}\notag\\ 
&\leq C_2\big(\|{\mathcal E}_0'(\omega 
+\eta)-{\mathcal E}_0''(\omega)\eta\|_{W^{1,3}}+2|\e|\|{\mathcal
  V}_H'(\omega+p+\eta)\|_{W^{1,3}}\big)\notag\\ 
&\leq C_2\bigg(\int_0^1\|{\mathcal E}_0''(\omega 
+t \eta)-{\mathcal E}_0''(\omega)\|_{W^{1,3/2}}\,dt\|\eta\|_{W^{1,3}}+
2|\e|\|{\mathcal V}_H'(\omega+p+\eta)\|_{W^{1,3}}\bigg)\notag\\
&\leq C_2\rho\sup_{\|\eta\|_{W^{1,3}}\leq\rho}\|{\mathcal E}_0''(\omega
+ \eta)-{\mathcal E}_0''(\omega)\|_{W^{1,3/2}}+2
C_2|\e|\sup_{\|\eta\|_{W^{1,3}}\leq\rho}\|{\mathcal
  V}_H'(\omega+p+\eta)\|_{W^{1,3}} 
\end{align}
where $C_2=\|{\mathcal L}^{-1}\|_{{\mathcal L}(W^{1,3}\times\R^6\times\R^3)}$. For $(\eta_1,\l_1,\a_1)$, $(\eta_2,\l_2,\a_2)\in B_{\rho}(0)\subset W^{1,3}\times\R^6\times\R^3$ we have
\begin{align*}
&\frac{\|T_{\e,p}(\eta_1,\l_1,\a_1)-T_{\e,p}(\eta_2,\l_2,\a_2)\|_{W^{1,3}\times\R^6\times\R^3}}{C_2 \|\eta_1-\eta_2\|_{W^{1,3}}}\\[5pt]
&\quad\leq \frac{\|{\mathcal E}_{\e}'(\omega+p+\eta_1)-{\mathcal E}_{\e}'(\omega+p+\eta_2)-{\mathcal E}_0''(\omega)(\eta_1-\eta_2)\|_{W^{1,3}}}{C_2 \|\eta_1-\eta_2\|_{W^{1,3}}}\\
&\quad\leq \int_0^1\|{\mathcal E}_{\e}''(\omega+p+\eta_2
+t (\eta_1-\eta_2))-{\mathcal E}_0''(\omega)\|_{W^{1,3/2}}\,dt\\
&\quad\leq \int_0^1\|{\mathcal E}_{0}''(\omega+p+\eta_2
+t (\eta_1-\eta_2))-{\mathcal E}_0''(\omega)\|_{W^{1,3/2}}\,dt\\
&\quad\quad+
2|\e|\int_0^1\|{\mathcal V}_H''(\omega+p+\eta_2+t(\eta_1-\eta_2))\|_{W^{1,3/2}}\,dt\\
&\quad\leq\sup_{\|\eta\|_{W^{1,3}}\leq 3\rho}\|{\mathcal E}_0''(\omega
+ \eta)-{\mathcal E}_0''(\omega)\|_{W^{1,3/2}}+2
|\e|\sup_{\|\eta\|_{W^{1,3}}\leq 3\rho}\|{\mathcal
  V}_H''(\omega+p+\eta)\|_{W^{1,3/2}}. 
\end{align*}
 From \eqref{eq:frechet1}, \eqref{eq:frechet2}, and the H\"older inequality it follows that there exists a positive constant $C_3$ such that for any $\eta\in W^{1,3}$, $p\in\R^3$
\begin{align}
 \|{\mathcal V}_H'(\omega+p+\eta)\|_{W^{1,3}}&\leq
 C_3\bigg[\bigg(\int_{\R^2}|H(\omega+p+\eta)|^{3/2}|\n\omega|^3\mu^{-1}\bigg)^{2/3}+\|\eta\|^2_{W^{1,3}}\bigg]\label{eq:est1}\\ 
\|{\mathcal V}_H''(\omega+p+\eta)\|_{W^{1,3/2}}&\leq
C_3\bigg[\bigg(\int_{\R^2}|H(\omega+p+\eta)|^2|\n(\omega+\eta)|^2\bigg)^{1/2}\notag\\
&\qquad\quad+ \bigg(\int_{\R^2}|\n
H(\omega+p+\eta)|^{3/2}|\n(\omega+\eta)|^3\mu^{-1}\bigg)^{2/3}
\bigg]\label{eq:est2}. 
\end{align}
Choosing $\rho_0>0$ such that 
\[
C_2\sup_{\|\eta\|_{W^{1,3}}\leq3\rho_0}\big\|{\mathcal E}_0''(\omega+\eta)-{\mathcal E}_0''(\omega)\big\|_{W^{1,3/2}}<\frac12
\]
and $\e_0>0$ such that
\begin{align}
8C_2C_3\e_0\|H\|_{L^{\infty}(\R^3)}\|\omega\|^2_{W^{1,3}}&< \min\bigg\{ 1,\rho_0,\frac1{8 C_2C_3\e_0}\bigg\},\label{eq:chose1}\\
\sup_{\substack{\|\eta\|_{W^{1,3}}\leq\rho_0\\p\in\R^3}}
\big\|{\mathcal V}_H'(\omega+p+\eta)\|_{W^{1,3}}& <\frac{\rho_0}{6\e_0 C_2},\label{eq:chose2}\\
\sup_{\substack{\|\eta\|_{W^{1,3}}\leq3\rho_0\\p\in\R^3}}\big\|{\mathcal V}_H''(\omega+p+\eta)\|_{W^{1,3/2}}& <\frac{1}{8\e_0 C_2},\label{eq;chose3}
\end{align}
we obtain that $T_{\e,p}$ maps the ball $\overline{B_{\rho_0}(0)}$ into itself
for any $|\e|< \e_0$, $p\in\R^3$, and is a contraction there. Hence it has a unique fixed point $\big(\eta(\e,p),\l(\e,p),\a(\e,p)\big)\in \overline{B_{\rho_0}(0)}$. From \eqref{eq:1} we have that the following property holds 
\begin{align*}
(*)\qquad&\text{$T_{\e,p}$ maps a ball $\overline{B_{\rho}(0)}\subset W^{1,3}\times \R^6\times \R^3$ into itself whenever $\rho\leq \rho_0$ and}\\
&\qquad\qquad
\rho> 4|\e|C_2\sup_{\|\eta\|_{W^{1,3}}\leq\rho}\big\|{\mathcal V}_H'(\omega+p+\eta)\|_{W^{1,3}}.
\end{align*}
In particular let us set 
\begin{equation}\label{eq:rhoeps}
\rho_{\e}=5|\e|C_2\sup_{\substack{\|\eta\|_{W^{1,3}}\leq\rho_0\\p\in\R^3}}\big\|{\mathcal V}_H'(\omega+p+\eta)\|_{W^{1,3}}.
\end{equation}
In view of \eqref{eq:chose2} and \eqref{eq:rhoeps}, we have that for any $|\e|<\e_0$ and for any $p\in\R^3$
\[
\rho_{\e}\leq\rho_0\quad\text{and}\quad \rho_{\e}> 4|\e|C_2\sup_{\|\eta\|_{W^{1,3}}\leq\rho_{\e}}\big\|{\mathcal V}_H'(\omega+p+\eta)\|_{W^{1,3}}
\]
so that, due to $(*)$, $T_{\e,p}$ maps $\overline{B_{\rho_{\e}}(0)}$
into itself. From the uniqueness of the fixed point we have that for
any $|\e|<\e_0$ and $p\in\R^3$ 
\begin{equation}\label{eq:bound}
\|(\eta(\e,p),\l(\e,p),\a(\e,p))\|_{W^{1,3}\times \R^6\times \R^3}\leq\rho_{\e}\leq C_1|\e|
\end{equation}
for some positive constant $C_1$ independent of $p$ and hence
$
\|\eta(\e,p)\|_{W^{1,3}}\leq\rho_{\e}\leq C_1|\e|
$
thus proving \eqref{eq:eta4}. Assume now \eqref{eq:limitH} and set for any $p\in\R^3$
\[
\rho_p=8 C_2C_3\e_0\bigg(\int_{\R^2}\sup_{|q-p|\leq 1+C_0}|H(q)|^{3/2}|\n\omega|^3\mu^{-1}\bigg)^{2/3}
\]
where $C_0$ is a positive constant such that $\|u\|_{L^{\infty}}\leq
C_0\|u\|_{W^{1,3}}$ for any $u\in W^{1,3}$. From~\eqref{eq:chose1} we have that 
\[
\rho_p<\min\bigg\{1,\rho_0,\frac1{8 C_2 C_3 \e_0}\bigg\}.
\]
Hence, due to \eqref{eq:est1}, we have that for $|\e|<\e_0$ and $\|\eta\|_{W^{1,3}}\leq\rho_p$
\begin{align*}
&4|\e|C_2\big\|{\mathcal V}_H'(\omega+p+\eta)\|_{W^{1,3}}\\
&\qquad
\leq 4\e_0 C_2C_3
\bigg(\int_{\R^2}\sup_{|q-p|\leq 1+C_0}|H(q)|^{3/2}|\n\omega|^3\mu^{-1}\bigg)^{2/3}+4\e_0C_2C_3\rho_p^2<\rho_p.
\end{align*}
From $(*)$ and the uniqueness of the fixed point, we deduce that $\|\eta(\e,p)\|_{W^{1,3}}\leq\rho_p$ for any $|\e|<\e_0$ and $p\in \R^3$. On the other hand, since $H$ vanishes at $\infty$, by the definition of $\rho_p$ we have that $\rho_p\to 0$ as $|p|\to\infty$, hence 
\[
\lim_{|p|\to\infty}\eta(\e,p)=0\quad\text{in }W^{1,3}\text{ uniformly for }|\e|<\e_0.
\]
The proof of Lemma \ref{l:eta} is now complete.
\end{pf}
\begin{rem}\label{r:eta}
The map $\eta$ given in Lemma \ref{l:eta} satisfies
\[
\langle {\mathcal
  E}_{\e}'(\omega+p+\eta(\e,p)),\varphi\rangle-\sum_{i=1}^6\l_i(\e,p)\int_{\R^2}\n \varphi\cdot\n\tau_i+\a(\e,p)\cdot\int_{{\mathbb S}^2}\varphi,\quad\forall \varphi\in W^{1,3/2}
\]
where $\big(\eta(\e,p),\l(\e,p),\a(\e,p)\big)\in
\overline{B_{\rho_{\e}}(0)}\subset W^{1,3}\times\R^6\times\R^3$ being
$\rho_{\e}$ given in (\ref{eq:rhoeps}), hence 
\begin{align*} 
\int_{\R^2}&\n(\omega+\eta(\e,p))\cdot\n\varphi+2\int_{\R^2}\varphi
\cdot(\omega+\eta(\e,p))_x\wedge(\omega+\eta(\e,p))_y\\
&\quad+ 2\e\int_{\R^2}H(\omega+p+\eta(\e,p))\varphi
\cdot(\omega+\eta(\e,p))_x\wedge(\omega+\eta(\e,p))_y\\
&=\sum_{i=1}^6\l_i(\e,p)\int_{\R^2}\n \varphi\cdot\n\tau_i-\a(\e,p)\cdot\int_{{\mathbb S}^2}\varphi,\quad\forall \varphi\in W^{1,3/2},
\end{align*}
i.e. $\eta(\e,p)$ satisfies the equation
\begin{align*}
\D \eta(\e,p)=F(\e,p)
\end{align*}
where
\begin{align*}
F(\e,p)=&2 (\omega+\eta(\e,p))_x\wedge(\omega+\eta(\e,p))_y-2
\omega_x\wedge \omega_y+\l(\e,p)\cdot\D
\tau-\a(\e,p)\mu^2\\
&+2\e
H(\omega+p+\eta(\e,p))(\omega+\eta(\e,p))_x\wedge(\omega+\eta(\e,p))_y\quad\mbox{in
  }\R^2. 
\end{align*}
Since $F(\e,p)\in L^{3/2}$ and, in view of (\ref{eq:eta4}) and (\ref{eq:bound}),
$F(\e,p)\to 0$ in $L^{3/2}$ as $\e\to0$ uniformly with respect to $p$, by regularity we have that 
\[
\eta(\e,p)\in W^{2,3/2}\quad\text{and}\quad \eta(\e,p)\to 0\text{ in
  }W^{2,3/2}
\]
hence, by Sobolev embeddings, $F(\e,p)\in L^{3}$ and $F(\e,p)\to 0$ in
$L^{3}$ as $\e\to0$ uniformly with respect to $p$. Again by
regularity 
\[
\eta(\e,p)\in W^{2,3}\quad\text{and}\quad \eta(\e,p)\to 0\text{ in
  }W^{2,3}
\]
hence $\eta(\e,p)\in C^{1,1/3}$ and 
\begin{equation}\label{eq:c1alpha}
\eta(\e,p)\to 0 \quad\text{in }C^{1,1/3}\text{ 
as }\e\to0\text{ uniformly with respect to }p. 
\end{equation}
\end{rem}
\par\noindent
For any $\e\in(-\e_0,\e_0)$, let us define the perturbed manifold
\[
Z_{\e}:=\big\{\omega+p+\eta(\e,p):\ p\in\R^3\big\}.
\]
From \cite{CM}, we have that $Z_{\e}$ is a natural constraint for ${\mathcal E}_{\e}$, namely any critical point $p\in\R^3$ of the functional 
\[
\Phi_{\e}:\R^3\to\R,\qquad   \Phi_{\e}(p)={\mathcal E}_{\e}\big(\omega+p+\eta(\e,p)\big)
\]
gives rise to a critical point $u_{\e}=\omega+p+\eta(\e,p)$ of ${\mathcal E}_{\e}$.
\begin{pro}\label{p:phieps}
Assume $H\in C^2(\R^3)$, $\n H\in  L^{\infty}(\R^3,\R^3)$, and
$\lim_{|p|\to\infty}H(p)=0$. Then for any $|\e|<\e_0$
\[
\lim_{|p|\to\infty}\Phi_{\e}(p)={\rm const}={\mathcal E}_0(\omega).
\]
\end{pro}
\begin{pf}
We have that
\begin{align}\label{eq:phieps1}
\Phi_{\e}(p)&={\mathcal E}_{\e}\big(\omega+p+\eta(\e,p)\big)={\mathcal E}_0
\big(\omega+p+\eta(\e,p)\big)+2\e{\mathcal V}_H\big(\omega+p+\eta(\e,p)\big)\notag\\
&=\frac12 \int_{\R^2}|\n \omega|^2+\frac12 \int_{\R^2}|\n \eta(\e,p)|^2+\int_{\R^2}\n\omega\cdot\n\eta(\e,p)\notag\\
&\quad+\frac23 
\int_{\R^2} \big(\omega+p+\eta(\e,p)\big)\cdot\big(\omega+\eta(\e,p)\big)_x\wedge \big(\omega+\eta(\e,p)\big)_y\notag\\
&\quad+2\e \big[{\mathcal V}_H(\omega+p)+\langle 
{\mathcal V}_H'(\omega+p),\eta(\e,p)\rangle+o\big(\|\eta(\e,p)\|_{W^{1,3}}\big)\big]\notag\\
&=\frac12 \int_{\R^2}|\n \omega|^2
+\frac23\int_{\R^2}\omega\cdot\omega_x\wedge\omega_y 
+\frac12 \int_{\R^2}|\n \eta(\e,p)|^2+\int_{\R^2}\n\omega\cdot\n\eta(\e,p)
\notag\\
&\quad+\frac23\int_{\R^2}\omega\cdot(\omega_x\wedge\eta(\e,p)_y+\eta(\e,p)_x\wedge\omega_y)\notag\\
&\quad
+\frac23\int_{\R^2}\omega\cdot\eta(\e,p)_x\wedge\eta(\e,p)_y  
+\frac23\int_{\R^2}\eta(\e,p)\cdot(\omega+\eta(\e,p))_x\wedge(\omega+\eta(\e,p))_y \notag\\
&\quad+2\e {\mathcal V}_H(\omega+p)+2\e  \langle 
{\mathcal V}_H'(\omega+p),\eta(\e,p)\rangle+2\e\, o\big(\|\eta(\e,p)\|_{W^{1,3}}\big)
\end{align}
where we have used the fact that
\[
\int_{\R^2}p\cdot u_x\wedge u_y=0\quad\forall\, p\in \R^3,\ u\in W^{1,3},
\]
(see \cite{CM}, Lemma A.3). Notice that from Lemma \ref{l:eta} we have that
\begin{align}
\int_{\R^2}|\n \eta(\e,p)|^2&\leq \sqrt[3]{4\pi}\|\eta(\e,p)\|_{W^{1,3}}^2\mathop{\longrightarrow}\limits_{|p|\to\infty} 0,\label{eq:phieps2}\\
\bigg|\int_{\R^2}\n \omega\cdot\n \eta(\e,p)\bigg|&\leq \sqrt[6]{4\pi}\bigg(\int_{\R^2}|\n\omega|^2\bigg)^{1/2}\|\eta(\e,p)\|_{W^{1,3}}\mathop{\longrightarrow}\limits_{|p|\to\infty} 0,\label{eq:phieps3}
\end{align}
and, by the H\"older inequality and Lemma \ref {l:eta},
\begin{align}
&\bigg|\int_{\R^2}\omega\cdot(\omega_x\wedge\eta(\e,p)_y+\eta(\e,p)_x\wedge\omega_y)\bigg|\leq 2\|\omega\| _{W^{1,3}}^2\|\eta(\e,p)\|_{W^{1,3}}\mathop{\longrightarrow}\limits_{|p|\to\infty} 0,\label{eq:phieps4}\\
&\bigg|\int_{\R^2}\omega\cdot(\eta(\e,p)_x\wedge\eta(\e,p)_y)\bigg|\leq \|\omega\| _{W^{1,3}}\|\eta(\e,p)\|_{W^{1,3}}^2\mathop{\longrightarrow}\limits_{|p|\to\infty} 0,\label{eq:phieps5}\\
&\bigg|\int_{\R^2}\eta(\e,p)\cdot(\omega+\eta(\e,p))_x\wedge(\omega+\eta(\e,p))_y \bigg|\notag\\
&\hskip5truecm\leq \|\omega+\eta(\e,p)\| _{W^{1,3}}^2\|\eta(\e,p)\|_{W^{1,3}}\mathop{\longrightarrow}\limits_{|p|\to\infty} 0.\label{eq:phieps6}
\end{align}
Moreover the Gauss-Green Theorem yields
\[
{\mathcal V}_H(\omega+p)=-\int_{B_1}H(\xi+p)\,d\xi
\]
so that by the Dominated Convergence Theorem we have that
\begin{equation}\label{eq:phieps7}
\lim_{|p|\to\infty}{\mathcal V}_H(\omega+p)=0.
\end{equation}
From \eqref{eq:frechet1}, H\"older inequality, and Lemma \ref{l:eta}, we have that 
\begin{align}\label{eq:phieps8}
|\langle {\mathcal V}_H'(\omega+p),\eta(\e,p)\rangle|&=\bigg|\int_{\R^2} H(\omega+p)\,\eta(\e,p)\cdot\omega_x\wedge\omega_y\bigg|\notag\\
&\leq \| H\|_{L^{\infty}(\R^3)}\|\omega\|^2_{W^{1,3}}\|\eta(\e,p)\|_{W^{1,3}}
\mathop{\longrightarrow}\limits_{|p|\to\infty} 0.
\end{align}
From \eqref{eq:phieps1} - \eqref{eq:phieps8}, it follows that
\[
\lim_{|p|\to\infty}\Phi_{\e}(p)=\frac12 \int_{\R^2}|\n \omega|^2+\frac23 
\int_{\R^2}\omega\cdot \omega_x\wedge\omega_y={\mathcal E}_0(\omega).
\]
The proposition is thereby proved.
\end{pf}

\begin{pfn}{Theorem \ref{t:1}.}
As already observed at the beginning of Section 2, it is not
restrictive to take $H_0=1$. From Proposition \ref{p:phieps} it
follows that for $|\e|<\e_0$ either $\Phi_{\e}$ is constant (and hence
we have infinitely many critical points) or it has a global maximum or
minimum point. In any case $\Phi_{\e}$ has a critical point. Since
$Z_{\e}$ is a natural constraint for ${\mathcal E}_{\e}$, we deduce
the existence of a critical point of ${\mathcal E}_{\e}$ for
$|\e|<\e_0$ and hence of a solution to $(P_{\e})$. The $H_{\e}$-bubble
$\omega_{\e}$ 
found in this way is of the form $\omega+p^{\e}+\eta(\e,p^{\e})$ for
some $p^{\e}\in\R^3$ where $\eta$ is as in Lemma \ref{l:eta}. 
Remark \ref{r:eta} yields  that $\omega_{\e}$ is closed
in~$C^{1,1/3}({\mathbb S}^2,\R^3)$-norm to the manifold $\{\omega+p:\ 
p\in \R^3\}$ for $\e$  
small. Since $\omega$ has no branch points, we deduce that $\omega_{\e}$ has no
branch points. 
\end{pfn}

\noindent To prove Theorems \ref{t:2} and \ref{t:3}, we need the following expansion for
$\Phi_{\e}$ (see \cite{CM})
\begin{equation}\label{eq:expansion}
\Phi_{\e}(p)={\mathcal E}_0(\omega)-2\e \Gamma(p)+O(\e^2)\quad\text{as
  }\e\to 0\text{ uniformly in }p\in\R^3.
\end{equation}

\begin{pfn}{Theorem \ref{t:2}.} Let $\e>0$ small. Assumption $(H4)$
  implies that 
  $\Gamma(0)>0$ and hence from (\ref{eq:expansion}) we have that for
  $\e$ small  $\Phi_{\e}(0)<{\mathcal E}_0(\omega)$, whereas from
  assumption $(H3)$ we have that $\hes\Gamma(0)$ is positive definite
  so that $\Gamma$ has a strict local minimum  
  in $0$ and hence from (\ref{eq:expansion}) $\Phi_{\e}$ has a strict
  local maximum in $B_r(0)$ for some $r>0$ such that
  $\Phi_{\e}(p)<\Phi_{\e}(0)-c_{\e}<{\mathcal E}_0(\omega)$ for $|p|=r$,
  where $c_{\e}$ is some 
  positive constant depending on $\e$. In particular $\Phi_{\e}$ has a
  mountain pass geometry. Moreover by Theorem \ref{t:1}
  $\Phi_{\e}(p)\to {\mathcal E}_0(\omega)$ as $|p|\to\infty$, and so
  $\Phi_{\e}$ must have a global minimum point. If the minimum point
  and the mountain pass point coincide then $\Phi_{\e}$ has infinitely
  many critical points. Otherwise $\Phi_{\e}$ has at least
  three critical points: a local maximum point, a global minimum
  point, and a mountain pass. If $\e<0$ we find the inverse
  inequalities and hence we find that $\Phi_{\e}$ has a local minimum
  point, a global maximum 
  point, and a mountain pass. As a consequence $(P_{\e})$ has at least
  three solutions provided $\e$ is sufficiently small.~\end{pfn}

As observed in Remark \ref{r:2}, if $H_1(0)>0$ and $\hes H_1(0)$ is
positive definite, by
continuity we
have that for $H_0$ sufficiently large $\Gamma(0)>0$ and
$\hes\Gamma(0)$ is positive definite, so that we can still prove the
existence of three
solutions arguing as above.

\begin{pfn}{Theorem \ref{t:3}.} Assumption $(H5)$ implies that
  $\Gamma(p_1)>0$ and $\Gamma(p_2)<0$. Since $\Phi_{\e}(p)={\mathcal
    E}_0(\omega)+2\e\big(- \Gamma(p)+o(1)\big)$ as $\e\to0$, we have
  for $\e$ sufficiently small
\[
\Phi_{\e}(p_1)<{\mathcal E}_0(\omega)\quad\text{and}\quad
\Phi_{\e}(p_2)>{\mathcal E}_0(\omega)
\]
if $\e>0$ and the inverse inequalities if $\e<0$.    
Since by Theorem \ref{t:1} $\Phi_{\e}(p)\to {\mathcal E}_0(\omega)$ as
$|p|\to\infty$, we conclude that $\Phi_{\e}$ must have a global maximum
point and a global minimum point in $\R^3$. Since $Z_{\e}$ is a
natural constraint for ${\mathcal E}_{\e}$, we deduce the existence of
two critical points of ${\mathcal E}_{\e}$ for $|\e|$ sufficiently
small and hence of two solutions to $(P_{\e})$. 
\end{pfn}

As observed in Remark \ref{r:3}, if $H_1(p_1)>0$ and $H_1(p_2)<0$, by
continuity we
have that for $H_0$ sufficiently large $\Gamma(p_1)>0$ and
$\Gamma(p_2)<0$, so that we can still prove the existence of two
solutions arguing as above.

\end{document}